\documentclass[10pt,twocolumn,twoside]{IEEEtran}

\usepackage{amsmath,amsfonts,mathrsfs}
\usepackage{graphicx}
\usepackage{theorem}

\usepackage{amssymb}

\usepackage{epsfig}

\usepackage{rotating}

\newcommand{\brk}[1]{\left( #1 \right)}
\newcommand{\Brk}[1]{\left[ #1 \right]}
\newcommand{\BRK}[1]{\left\{ #1 \right\}}

\renewcommand{\Re}{\mathbb{R}}

\newcommand{\e}{\mathbf{E}}

\newcommand{\eps}{\epsilon}
\newcommand{\inveps}{\frac{1}{\eps}}

\newcommand{\Dt}{\Delta t}
\newcommand{\Ds}{\Delta s}
\newcommand{\dt}{\delta t}
\newcommand{\dw}{\Delta W}

\newcommand{\bQ}{\bar Q}

\newcommand{\ba}{\bar{a}}

\newtheorem {algorithm}{Algorithm}[section]

\numberwithin{equation}{section}

\begin{document}

\title{Variance Reduction for Particle Filters of Systems with Time Scale Separation}
\author{Dror Givon, Panagiotis Stinis, and Jonathan Weare
\thanks{D. Givon, P. Stinis, and J. Weare are affiliated 
with the Department of Mathematics,
University of California, and
Lawrence Berkeley National Laboratory
    Berkeley, CA 94720}
\thanks{givon@math.berkeley.edu; stinis@math.lbl.gov; weare@cims.nyu.edu }
}

\maketitle

\begin{abstract}
We present a particle filter construction for a system that exhibits
time scale separation. The separation of time scales allows
two simplifications that we exploit: i) the use of the averaging
principle for the dimensional reduction of the dynamics for each
particle during the prediction step and ii) the factorization of the
transition probability for the Rao-Blackwellization of the update
step. The resulting particle filter is faster and
has smaller variance than the particle filter based on the original
system.  The method is tested on a multiscale stochastic differential equation
and on a multiscale pure jump diffusion motivated by chemical reactions.
\end{abstract}
\IEEEpeerreviewmaketitle


\section{Introduction}
The field of dimensional reduction has seen a flourishing in the last decade
(see e.g. the reviews in \cite{CS05,GKS04}), mainly due to i) the realization
that many systems of physical interest are more complex than one can handle
even with the largest available computers and ii) the fact that for many
complex systems the quantities of interest are coarse-scale features. Once a reduced model is constructed, it can be used in conjunction with filtering algorithms, like particle filters \cite{DFG01}, to incorporate information from real-time measurements. If the system under consideration exhibits time scale separation, the construction of a reduced model and subsequently of a particle filter, is also simplified. In this work we present a particle filter which exploits this simplification to create a particle
filter that is more efficient than the particle filter of the original
(large dimensional) system. The particle filter we construct is proved to
converge to the analytical filter of the original system.

A strategy for reducing the number
of unknowns in  an
application of a  filter for
nonlinear dynamical systems is proposed in \cite{CK04}. In that work the
authors assumed that the dynamics are strongly locally contractive
in some directions which are found adaptively.
Here we make an alternative assumption.  We instead assume that the
system exhibits a time scale separation. By this we mean that
certain components or modes of the system tend to move very slowly
in comparison with the rest of the modes.
In particular we are
interested in approximating conditional expectations the form
$$
\mathbf{E}\left[f\left(X_{s_k}\right)
\vert \left\{ Z_1,\dots,Z_k \right\}\right]
$$
where
$
Z_k$
are noisy observations of some multiscale Markov process $X_t$ at discrete times
$s_1,\dots,s_N.$
For a very general definition of a multiscale Markov process see \cite{kur73,Pap77} and \cite{ps08}.
We will illustrate the method presented below for the particular case that $X_t$ is the solution to 
a stochastic differential equations with time scale
separation (see equation \eqref{eq:system}).
The basic idea, however, can be applied to filtering problems for a variety of multiscale Markov process for which fast multiscale integration methods are available (see \cite{GK03,Van03,ELV05,GKK06,CGP05,elv07,givon:495,ps08} and the 
references therein).  In section \ref{sec:numbers}, our filtering method
is applied first to the reconstruction of the trajectory of a multiscale stochastic differential equation and second to the reconstruction of a trajectory of a pure jump Markov process motivated by chemical reactions that takes place on vastly different time scales.

Loosely speaking, recursive estimation of conditional expectations of the form above
can be accomplished in two steps.  First, in the prediction step,
the system is evolved according to its evolution law. Second, in the update step, the resulting samples of
the system are weighted by the likelihood of the next observation
given the sample. We describe the filtering
problem in more detail in the next section.
The separation of
time scales facilitates the construction of an efficient particle
filter in two ways: i) it allows for fast evolution of the system in
the prediction step and ii) it allows the integration of the observation
weights over the fast modes of the system during the update step. Step ii)
amounts to a Rao-Blackwellization of a standard particle filter estimator.

Multiscale phenomena have been observed in wide ranging areas of research.  For example,
empirical evidence from a study of exchange rate dynamics in \cite{abd02},
suggests the use of stochastic volatility models with both fast and slow time scales.
Other examples of systems with scale separation include
chemical reaction systems where there can exist a difference of several orders of magnitude among the different reaction rates (see \cite{CGP05,elv07,hr02,ra03,rmgk04}). 
Similar problems exist in material science (see \cite{curtin}) and molecular dynamics simulations (see \cite{mttk96}) where one is interested in large scale features of the system but this behavior depends critically on the small (and fast) scale motion. An even more challenging problem is that of weather prediction and how to assimilate (through filtering) the vast amount of measurements collected daily around the globe. The weather system exhibits an extremely large range of active scales not necessarily with clear time scale separation (see \cite{miller,weare08}). A projection formalism
framework, for the construction of reduced models of large systems with or
without time scale separation, has been presented by Chorin and co-workers
(see \cite{chorin1}).

The main difficulty presented by multiscale phenomena is that they are extremely costly to integrate.  The vast majority of computational time is spent evolving the fast components of the system while one may be primarily interested in slow scale characteristics.
This implies that the prediction step in any filtering method which does not take directly into account the multiscale structure of the
problem will become prohibitively costly as the time scale separation increases.  As described in detail below, our method
addresses this issue through its incorporation of the multiscale
integration scheme.

A second issue, common to importance sampling techniques for high dimensional problems, 
is controlling the variance of the resulting estimator.  To address this in the context of particle filtering, 
many authors have suggested the use of some form of Rao-Blackwellization
(see \cite{ah77,dga00,dgk01,klw94,lc98,Vas07} for example).
The distribution of the underlying Markov process given current and past observations 
can always be factored into a posterior marginal distribution for
some set of (in our case slow) variables and a posterior conditional distribution for the remaining
(in our case fast) variables given the first set of variables.
Rao-Blackwellization requires that expectations with respect to 
the resulting posterior conditional distribution
can be carried out exactly.  In the case that the posterior conditional distribution
can be sampled, the Rao-Blackwellization procedure
 can be  approximated
by Monte Carlo averaging over these samples. 
As discussed in detail later, the separation of scales assumption made in
 this paper allows an approximate factorization of the posterior distribution in which
samples from the posterior conditional distribution can be easily and efficiently generated.  
While our assumption on the system clearly restricts
the class of possible applications, another advantage of our setup is
that the posterior conditional distribution of the unresolved variables is
allowed to be quite general (for example very non-Gaussian).


A method closely related to the one presented here was recently
proposed in \cite{Pap08}. There is, however, an important difference. 
In \cite{Pap08} it is assumed
that the observations and the objective function ($f$ in the conditional expectation above) depend directly only on the slow modes in the
system.  Here we allow for general observations and general
objective functions.  This fact is central to
the utility of our algorithm.  The method proposed here is designed to not
only improve the efficiency of the prediction step of a particle filter, but also 
to reduce the number of particles required to achieve a given accuracy.  In fact, in some problems with a moderate time scale separation, one may not observe any gain in efficiency in the
prediction step, while the reduction in variance may be significant.
It should also be noted that analytical results for continuous time multiscale filtering problems have been obtained by Kushner in
\cite{Kus90}.

The paper is organized as follows. Section
\ref{sec:filteringproblem} recalls the well known particle filter
construction for a general Markov process which is observed with noise
 at a discrete set of instants. Section
\ref{sec:reduced} presents our particle filter construction in the
particular case of a multiscale stochastic differential equation.
Section \ref{sec:scalerao} discusses how the presence of a
separation of time scales can lead to a construction of a reduced
model for the slow components of the system and how it allows the
Rao-Blackwellization of the filtering step. It also includes a
particle filter construction which is based on the assumption that
the reduced model can be constructed analytically and the
Rao-Blackwellization of the filtering step can be performed
analytically. Section \ref{sec:free} contains the main algorithm,
which approximates the particle filter presented in Section
\ref{sec:scalerao} when the necessary calculations cannot be
performed analytically. Section
\ref{sec:sde} contains numerical results for a system of 
stochastic differential equations.
Section
\ref{sec:jump} contains numerical results for a pure jump type Markov process
motivated by multiscale chemical reactions.
Finally, Section \ref{sec:conclusions} contains a discussion of the
algorithm and of the results.


\section{The filtering problem and particle filters}\label{sec:filteringproblem}
We start by formulating the filtering problem. Assume that $X_t$ is a $d$-dimensional Markov process with transition probability $Q_{\Ds}(x,dx')$ where $Q_0(dx)$ is a known distribution on $\Re^d$.  We assume for notational convenience, that the process is autonomous.  The process $X_t$ is usually called the hidden signal. Assume, also, that we have noisy discrete observations at $N$ regularly spaced times of length $\Delta s=s_k-s_{k-1}$ for $k=1,\ldots,N$ where $s_0=0$ and $s_k$ the time of the $i$-th observation. The observations satisfy,
\[
Z_k=G(X_{s_k},\chi_k), \quad k=1,\ldots,N
\]
where the $\chi_k$ are i.i.d random variables, independent of $X_t$, and for all $x$ the variable $Z=G(x,\chi)$ admits a density $z\to g(x,z)$ which is known. Let $1:k$ denote the sequence $\BRK{1,\ldots,k}$ for $k=1,\ldots,N$. The filtering problem  consists of computing the conditional expectations
\[
\Pi_k f = \e\Brk{f(X_{s_k}) | Z_{1:k}} \; \text{for} \; k=1,\ldots,N,
\]
where $f$ belongs to some reasonable family of functions. The quantities $\Pi_k f$ constitute the filter (called the analytical filter).

The conditional expectation $\Pi_k f$ can be written as
\begin{multline}\label{eq:filter1}
\Pi_k f =
\frac{\int Q_0(dx_0) Q_{\Delta s}(x_0,dx_{s_1})g(x_{s_1},z_1)}
{\int Q_0(dx_0) Q_{\Delta s}(x_0,dx_{s_1})g(x_{s_1},z_1)}\\
\cdots
\frac{Q_{\Delta s}(x_{s_{k-1}},dx_{s_k})g(x_{s_k},z_k)f(x_{s_k})}
{Q_{\Delta s}(x_{s_{k-1}},dx_{s_k})g(x_{s_k},z_k)},
\end{multline}
where $\int$ stands for integration over the variables $x_{s_0},x_{s_1},\dots,x_{s_k}.$
Let $H_k$ be the kernels,
\begin{multline}\label{eq:kernel}
H_k f\,(x_{s_{k-1}},z_k) = \int Q_{\Delta s}(x_{s_{k-1}},dx_{s_k})\\ \times g(x_{s_k},z_k)\ f(x_{s_k})\ dx_{s_k}.
\end{multline}
The filter (\ref{eq:filter1}) can be written recursively as
\begin{equation}\label{eq:filter2}
\Pi_0 f=\int f(x)\ Q_0(dx),\quad
\Pi_kf=\frac{\Pi_{k-1}H_k\,f}{\Pi_{k-1}H_k\,1}.
\end{equation}
Usually, the integrals in (\ref{eq:filter2}) cannot be computed
analytically. A particle filter \cite{GSS93, DFG01} is an
approximation to the analytical filter (\ref{eq:filter2}).  In its
simplest form, a particle filter consists of the following steps:

\begin{algorithm}{Standard particle filter.}\label{pf}

\begin{enumerate}
\item Sample $n$ i.i.d. vectors (particles) $X_0^j, \; j=1,\ldots,n$ from $Q_0(dx).$ Set $k=1.$
\item For $j=1,\ldots,n$ evolve $X_{s_{k-1}}^j$ according to $Q_{\Delta s}(X_{s_{k-1}}^j,dx')$ to get ${X'}_{s_k}^j$.
\item Calculate
    \[
    \begin{split}
    W^n_k\,f &= \frac{1}{n}\sum_{j=1}^n f({X'}_{s_k}^j)\, g({X'}_{s_k}^j,z_k), \\
    W^n_k\,1&= \frac{1}{n}\sum_{j=1}^n \, g({X'}_{s_k}^j,z_k).
    \end{split}
    \]
\item Define the probability measure $\Psi^n_k$ on $\Re^d$ with density
    \[
    \Psi^n_k(dx) = \frac{1}{n\,W^n_k\,1}\sum_{j=1}^n\ g({X'}_{s_k}^j,z_k)\ \delta_{{X'}_{s_k}^j}(dx).
    \]
\item If $k < N$: Set $ k \rightarrow k+1.$  Sample $n$ i.i.d. variables ${X}_{s_k}^j, \; j=1,\ldots,n$ from $\Psi^n_k$, and  go to Step 2.
\end{enumerate}
We stop at the end of iteration $N$, and our approximation of $\Pi_k\,f$ is,
\[
\Psi_k^n\,f = \frac{W_k^n\,f}{W_k^n\,1} = \e_{\Psi_k^n}f
\]
\end{algorithm}




\section{Filtering for systems with time scale separation}\label{sec:reduced}

\subsection{Scale separation and  Rao-Blackwellization} \label{sec:scalerao}
Many problems in the natural sciences give rise to systems with
time scale separation. These systems represent a challenge for
numerical simulations. For example in molecular dynamics simulations
femtoseconds timesteps are required to integrate the fastest atomic
motions, while the object of interest is of orders of microseconds
or longer. In the past four decades systems with scale separation
have been the focus of extensive research within the framework of
the averaging principle. The separation of scales is utilized to
derive an effective model for the slow components of the system.
While the averaging principle and its resulting effective dynamics
provide a substantial simplification of the original system it is
often impossible or impractical to obtain the reduced equations in
closed form. This has motivated the development of algorithms such
as multiscale integration methods described in the next section.

To make the presentation more concrete we will describe the situation for a $d$-dimensional systems of stochastic differential equations with multiple time scales (see \cite{FW84} Chapter 7).  As mentioned in the introduction, the ideas that we will 
discuss can be applied to filtering problems for any  multiscale Markov process for which multiscale integration methods are available.

Let $(X^{\epsilon}_t,Y^{\epsilon}_t)$ be a solution of the
system
\begin{subequations}
\label{eq:system}
\begin{align}
dX^{\epsilon}_t  &= a(X^{\epsilon}_t,Y^{\epsilon}_t)\,dt + b(X^{\epsilon}_t)\,dU_t,
 &X^{\epsilon}_0&=x_0 \label{eq:system1} \\
dY^{\epsilon}_t &= \inveps \alpha(X^{\epsilon}_t,Y^{\epsilon}_t)\,dt +
\frac{1}{\sqrt{\eps}}\beta(X^{\epsilon}_t,Y^{\epsilon}_t)\,dV_t,
 &Y^{\epsilon}_0&=y_0  \label{eq:system2}.
\end{align}
\end{subequations}
where $U_t$ and $V_t$ are independent Brownian motions.
The hidden variable $(X^{\eps}_t,Y^{\eps}_t)$ is a Markov process with transition probability $Q^{\eps}_{\Ds}((x,y),(dx',dy'))$ on 
$\Re^{d_x} \times \Re^{d_y}$ (with $d_x+d_y=d$). The $X^{\eps}_t$ variables evolve  on a $O(1)$ time scale (the macro time scale), and the $Y^{\eps}_t$  variables  evolve on an $O(\eps)$ time scale (the micro time scale). As above, we have noisy discrete observations,
\[
Z_k=G(X^{\eps}_{s_k},Y^{\eps}_{s_k},\chi_k), \; k=1,\ldots,N
\]
where the $\chi_k$ are i.i.d variables, independent of $x,y$, and for all $x,y$ the variable $Z=G(x,y,\chi)$ admits a density $z\to g(x,y,z)$ which is known.

The standard approach to this filtering problem (see \cite{DJP01}) is to use an easily sampled
approximation $Q^{\eps,\dt}_{\Ds}\brk{\brk{x,y},\brk{dx',dy'}}$ of the
transition probability $Q^{\eps}_{\Ds}((x,y),(dx',dy'))$ 
where the discrete time step $\dt$ is of scale comparable to $\eps.$ If one is interested in
the evolution of the system over $O(1)$ time scales then one must
evolve the system for $O(\frac{1}{\epsilon})$ steps, which can be
very costly for systems with large scale separation. In addition to simulation issues inherent to multiscale phenomena, particle filters can suffer from all of the usual difficulties associated with importance sampling in large dimensional spaces.
That is, in high dimensional systems, the variance of the particle filter estimator can be difficult to control.

 In this section we
show how the averaging principle and Rao-Blackwellization  can be
used to reduce the
computational effort for each particle and the number of required particles.  
The method we discuss
assumes that the Rao-Blackwellization and the construction of
the dimensionally reduced model can be performed analytically.
Unfortunately, this is rarely the case. In the next section, we show
how the multiscale integration framework can be used to implement
our approach when we are not able to perform the
Rao-Blackwellization and the construction of the reduced model
analytically.

Under appropriate assumptions (see \cite{Kha68,kur73,Pap77,FW84}), the averaging principle
dictates that as $\eps\to0$
\[
X^{\epsilon}_t \overset{\mathscr{D}}{\longrightarrow} \; \bar{X}_t \;  \text{for} \;  t\in[0,T],
\]
where $\bar{X}_t$ satisfies
\begin{equation}
\label{eq:effective}
d\bar{X}_t = \ba(\bar{X}_t)\,dt + b(\bar{X}_t)\,dW_t.
\end{equation}
The averaged coefficient $\ba$ is given by
\begin{equation}
\label{eq:averaging}
\begin{aligned}
\ba(x)&= \int_{\Re^{d_y}} a(x,y)\mu_x (dy)\\
\end{aligned}
\end{equation}
where $\mu_x (dy)$ is the invariant measure induced by \eqref{eq:system2} with the $x$ variables fixed.

The key consequence of the time scale separation is that, loosely speaking,
for $\eps$ small enough and for $\Ds\gg\eps,$ the transition density
$Q^{\eps}_{\Ds}\brk{(x,y),(dx',dy')}$ can be approximately factored as
\begin{equation}\label{eq:transapprox}
Q^{\eps}_{\Ds}\brk{(x,y),(dx',dy')} \approx \bQ_{\Ds}(x,dx')\,\mu_{x'}(dy)
\end{equation}
where $\bQ_{\Ds}$ is the transition density for the averaged system, \eqref{eq:effective}.
The factorization \eqref{eq:transapprox} above is the central tool in the construction of the multiscale particle filter below and is not a feature limited to multiscale stochastic differential equations.  The algorithms
below apply to any problem for which this approximation holds and for
which some means of  sampling (or approximately sampling)  $\bQ_{\Ds}(x,dx')$ and $\mu_{x'}(dy)$ are available.  This is the case for both of the examples in Section \ref{sec:numbers}.

The relation \eqref{eq:transapprox} suggests that an approximate particle filter can be constructed by defining the kernel,
\begin{multline}\label{eq:kernelaveraged}
\bar H_k f\,(x_{s_{k-1}},z_k) = \int \bar Q_{\Delta s}(x_{s_{k-1}},dx_{s_k})\\ \times\mu_{x_{s_{k}}}(dy)\ g(x_{s_k},y,z_k)\ f(x_{s_k},y).
\end{multline}
The corresponding particle filter would proceed as follows,
\begin{algorithm}{Standard particle filter corresponding to \eqref{eq:kernelaveraged}.}\label{avefilter}

\begin{enumerate}
\item Sample $n$ i.i.d. vectors $( X_0^j,Y_0^j) \; j=1,\ldots,n$ from $Q_0(dx,dy)$ distribution. Set $k=1.$

\item For each $j,$ use   \eqref{eq:averaging} to evolve
$X_{s_{k-1}}^j$ to ${X'_{s_k}}^j,$ i.e., sample from
$\bQ_{\Delta s}(X_{s_{k-1}}^j,dx').$  For each sample ${X'_{s_k}}^j,$ generate an independent sample
${Y'}^j$ from the measure $\mu_{{X'_{s_k}}^j}.$

\item   \label{item:loop}  

   Calculate,
    \[
    \begin{split}
    W^n_k\,f &= \frac{1}{n}\sum_{j=1}^n f ({X'_{s_k}}^j)\,g({X'_{s_k}}^j,{Y'}^j,z_k)\\
    W^n_k\,1 &= \frac{1}{n}\sum_{j=1}^n g({X'_{s_k}}^j,{Y'}^j,z_k).
\end{split}
    \]

\item   Define the probability measure $\Psi^n_k$ on $\Re^{d}$ with density
    \begin{multline*}
    \Psi^n_k(dx,dy) = \frac{1}{n W^n_k 1}\\\times \sum_{j=1}^n\
    g({X'_{s_k}}^j,{Y'}^j,z_k)\ \delta_{{X'_{s_k}}^j}(dx)\ \delta_{{Y'}^j}(dy) .
    \end{multline*}

\item If $k < N$: Set $ k \rightarrow k+1.$  Sample $n$ i.i.d. variables $X_{s_k}^j, \; j=1,\ldots,n$ from $\Psi^n_k$, and  go to Step 2.
\end{enumerate}
We stop at the end of iteration $N$, and our approximation of $\Pi_k\,f$ is,
\[
\Psi_k^n\,f = \frac{W_k^n\,f}{W_k^n\,1} = \e_{\Psi_k^n}f
\]
\end{algorithm}
There is no reason to hope that this algorithm should provide any variance reduction.  Indeed, in the ideal case that expression \eqref{eq:transapprox}
is an equality, the variance of the particle filter would be identical to
the variance of the standard particle filter \ref{eq:filter2}.

However, knowledge of  $\mu_{x_{s_{k}}}(dy)$ allows one to integrate out $y$ in expression \eqref{eq:kernelaveraged}  and write this same kernel in another form,
\begin{multline}\label{eq:kernelaveragedRB}
\bar H_k f\,(x_{s_{k-1}},z_k) = \int \bar Q_{\Delta s}(x_{s_{k-1}},dx_{s_k})\\ \times\int\mu_{x_{s_{k}}}(dy)\ g(x_{s_k},y,z_k)\ f(x_{s_k},y).
\end{multline}
This suggests a different filtering algorithm:
\begin{algorithm}{Rao-Blackwellized particle filter corresponding to \eqref{eq:kernelaveragedRB}.}
\label{avefilterrb}

\begin{enumerate}
\item Sample $n$ i.i.d. vectors $( X_0^j,Y_0^j) \; j=1,\ldots,n$ from $Q_0(dx,dy)$ distribution. Set $k=1.$

\item For each $j,$ use   \eqref{eq:averaging} to evolve
$X_{s_{k-1}}^j$ to ${X'_{s_k}}^j,$ i.e., sample from
$\bQ_{\Delta s}(X_{s_{k-1}}^j,dx').$

\item   \label{item:loop} Calculate the following quantities:   
    \begin{enumerate}
    \item\label{item:RB} 
     \begin{align*}
    \left[\mu fg\right]({X'_{s_k}}^j)& = \int f({X'_{s_k}}^j,y)\,
	\\ 
	&\hspace{0.5cm}\times g({X'_{s_k}}^j,y,z_k)\mu_{{X'_{s_k}}^j}(dy)\\
        \left[\mu g\right]({X'_{s_k}}^j) &= \int  g       ({X'_{s_k}}^j,y,z_k)\,\mu_{{X'_{s_k}}^j}(dy).
    \end{align*}
    \item   \label{item:RBreduced}
    \[
    \begin{split}
    W^n_k\,f &= \frac{1}{n}\sum_{j=1}^n \left[\mu fg\right]({X'_{s_k}}^j)\\
    W^n_k\,1 &= \frac{1}{n}\sum_{j=1}^n \left[\mu g\right]({X'_{s_k}}^j).
\end{split}
    \]
    \end{enumerate}

\item   Define the probability measure $\Psi^n_k$ on $\Re^{d}$ with density
    \begin{multline*}
    \bar{\Psi}^n_k(dx,dy) = \frac{1}{n\,W^n_k\,1}\\
	\times \sum_{j=1}^n\
    g({X'_{s_k}}^j,y,z_k)\ \delta_{{X'_{s_k}}^j}(dx)\ \mu_{{X'_{s_k}}^j}(dy) .
    \end{multline*}

\item If $k < N$: Set $ k \rightarrow k+1.$  Sample $n$ i.i.d. variables $X_{s_k}^j, \; j=1,\ldots,n$ from $\bar\Psi^n_k$, and  go to Step 2.
\end{enumerate}
We stop at the end of iteration $N$, and our approximation of $\Pi_k\,f$ is,
\begin{equation}\label{alg2key1}
\bar{\Psi}_k^n\,f = \frac{W_k^n\,f}{W_k^n\,1} = \e_{\bar\Psi_k^n}f
\end{equation}
\end{algorithm}

Since the $y$ components of the resampled particles after Step 4 are not used in Step 2, we
can, in practice, resample from the marginal density
\begin{equation}
\label{altresample1}
\bar{\Psi}^n_k(dx) = \frac{1}{n\,W^n_k\,1}\sum_{j=1}^n \
     \left[\mu g\right]\left({X'_{s_k}}^j\right)\ \delta_{{X'_{s_k}}^j}(dx).
\end{equation}

In order to illustrate the variance reduction aspects of algorithm \ref{avefilterrb} we will consider the asymptotic variance of 
a pair of related but much simpler estimators.  Define the conditional
expectations $\bar \Pi_k$ recursively by,
\begin{equation*}
\bar\Pi_0 f=\int f(x)\ Q_0(dx),\quad
\bar\Pi_kf=\frac{\bar\Pi_{k-1}\bar H_k\,f}{\bar \Pi_{k-1}\bar H_k\,1}.
\end{equation*}
Let
$$
I^n_k f = \frac{\frac{1}{n}\sum_{j=1}^n 
f(X^j_k,Y^j_k)\,g(X^j_k,Y^j_k,z_k)}
{\frac{1}{n}\sum_{j=1}^n g(X^j_k,Y^j_k,z_k)}
$$
and
$$
\bar I^n_k f = \frac{\frac{1}{n}\sum_{j=1}^n \left[\mu fg\right](X^j_k)}
{\frac{1}{n}\sum_{j=1}^n \left[g\right](X^j_k)}
$$
where $(X^{j}_k,Y^{j}_k)$ are $i.i.d.$ samples from the measure 
$\Pi_{k-1}\bQ_{\Ds}\mu$
Notice that the only difference between $I^n_k f$ and $\Psi^n_k f$ described
in Algorithm \eqref{eq:kernelaveraged} is that
the samples $(X^j_k,Y^j_k)$ are independently drawn from 
$\bar\Pi_{k-1}\bQ_{\Ds}\mu$
instead of its particle approximation $\Psi^n_{k-1}\bQ_{\Ds}\mu.$  
A corresponding relationship
holds between $\bar I^n_k f$ and $\bar\Psi^n_k f.$

The estimators $I^n_k f$ and $\bar I^n_k f$ satisfy Central Limit Theorems, i.e. 
\begin{equation}
 \sqrt{n} \left( I^n_k f - \bar\Pi_k f\right) 
\xrightarrow{\mathscr{D}} \mathcal{N}\left(0,\sigma^2\right)
\end{equation}
and
\begin{equation}
 \sqrt{n} \left( \bar I^n_k f - \bar \Pi_k f\right) 
\xrightarrow{\mathscr{D}} \mathcal{N}\left(0,\bar\sigma^2\right).
\end{equation}
An application of the delta method (see \cite{robert04}) yields that
\begin{equation*}
\sigma^2 = \bar\Pi_{k-1} \bQ_{\Ds}\mu\\
\left[ \left( f(X,Y) - \bar\Pi_k f\right)^2
\left(\frac{g(X,Y)}{\bar\Pi_{k-1} \bar H_k 1}\right)^2\right]
\end{equation*}
and
\begin{equation*}
\bar\sigma^2 = \bar\Pi_{k-1} \bar Q_{\Ds} 
\left[ \left( \frac{\left[\mu fg\right]({X})}{\left[\mu g\right](X)} -  \bar\Pi_k f\right)^2
\left(\frac{\left[\mu g\right](X)}
{\bar\Pi_{k-1}\bar H_k 1}\right)^2\right].
\end{equation*}
$\bar \sigma^2$ can be rewritten as
$$
 \bar\Pi_{k-1} \bar Q_{\Ds} 
\left[ \left[\mu \left(\frac{fg - g \bar\Pi_k f}{\bar\Pi_{k-1}\bar H_k 1}\right)\right]^2(X)\right].
$$
Therefore, by Jensen's inequality, we have that
$$
\bar \sigma^2 \leq 
 \bar\Pi_{k-1} \bar Q_{\Ds} 
\left[\mu \left(\frac{fg - g \bar\Pi_k f}{\bar\Pi_{k-1}\bar H_k 1}\right)^2\right] = \sigma^2.
$$

It is important to note that the variance reduction offered by 
Algorithm \ref{avefilterrb} does not require any Gaussian
or degenerate measure approximations.  The main assumptions are that
the system have a multiscale structure  and that one can 
sample
from $\bQ_{\Ds}(x,dx')$ and evaluate averages with respect to the ergodic measure $\mu_x$ (which we assume exists).  In many cases, for example
those in Section \ref{sec:numbers}, the first assumption can be
easily verified.
For most general systems, $\bQ_{\Ds}(x,dx')$ is not available in closed form.
In practice, we must replace $\bQ_{\Ds}$ by some approximation.  For example in the case of $\bar{x}_t$ above we can use the transition probability kernel 
for the Euler-Maruyama scheme with step size $\Dt,$ $\bQ^{\Dt}_{\Ds}.$
Of course to apply the Euler approximation to \eqref{eq:effective} we must 
be able to exactly evaluate averages with respect to $\mu_x.$  The removal of this assumption
will be addressed in the next section.


\subsection{Implementation of the multiscale integration for the reduced particle filter}\label{sec:free}

In the particle filter construction of the previous section  we used
the fact that we can average over the invariant measure induced by
the fast variables. This is usually impossible since  the invariant
measure is unknown or because integration cannot be performed
analytically. We will demonstrate that this problem can be overcome by
using multiscale
integration schemes (see \cite{GK03,Van03,ELV05,GKK06,CGP05,elv07,givon:495}).  
In the following description we will be discussing specifically the multiscale integration schemes of the form analyzed  in \cite{ELV05} and \cite{GKK06}.  The system studied
in section \ref{sec:jump} is a pure jump Markov process and therefore
requires a different multiscale integration scheme (see \cite{CGP05,elv07}).
Let
$\Dt$ be a fixed time step, and $X_{k,l}$ be the numerical
approximation to the coarse variable, $\bar{X}$ from the previous
section, at time $t_{k,l} = s_k+l\Dt$ (recall $s_k$ is the $k$-th
observation time). Assume for simplicity that $L=\frac{\Ds}{\Dt}$ is
an integer. Inspired by the limiting equation \eqref{eq:effective},
$X_{k,l}$ is evolved in time by an Euler-Maruyama step,
\begin{equation}
\label{eq:macro}
X_{k,l+1} = X_{k,l} + A(X_{k,l}) \,\Dt + b(X_{k,l}) \,\dw_{t_{k,l}}
\end{equation}
for $l=0,\dots,L-1,$
where $\dw_{t_{k,l}}$ are Brownian displacements over a time interval $\Dt$. We refer to \eqref{eq:macro} as the \emph{macro-solver}, or macro integrator, and  we denote its transition probability by $\bQ^{\Dt,\dt,M}_{\Ds}(x,dx').$  Notice that with this notation $X_{k,0} = X_{k-1,L}.$

The function $A(x)$ approximates $\ba(x),$ introduced in the previous section, 
which is the  
average of $a$ over an ergodic
measure. The ergodic property implies that
instead of ensemble averaging we can use time averaging over trajectories of
the rapid variables with fixed  $x$.
Since, by assumption, this average cannot be performed
analytically, it is approximated by an empirical average over
short runs of the fast dynamics. These ``short runs" are over time
intervals that are sufficiently long for empirical averages to be
close to their limiting ensemble averages, yet
sufficiently short for the entire procedure to be efficient
compared to the direct solution of the coupled system.

Thus, given the coarse variable at the ${k,l}$-th time step, $X_{k,l}$, we take some initial value for the fast component $Y_{k,l,0}$,
 and solve  \eqref{eq:system2} numerically with step size $\dt$ and $x=X_{k,l}$ fixed.
 We denote the discrete variables associated with the fast dynamics at the ${k,l}$-th coarse step by $Y_{k,l,m}$, $m=0,1,\dots,M.$ 
The numerical solver used to generate the sequence $Y_{k,l,m}$ is called the \emph{micro-solver}, or micro-integrator. The simplest choice is again the Euler-Maruyama scheme,
\begin{multline}
\label{eq:microsolver}
Y_{k,l,m+1} = Y_{k,l,m} + \inveps
\alpha(X_{k,l},Y_{k,l,m}) \,\dt\\ + \frac{1}{\sqrt{\eps}} \beta(X_{k,l},Y_{k,l,m}) \Delta
V_{k,l,m},
\end{multline}
where $\Delta V_{k,l,m}$ are Brownian displacements over a time
interval $\dt$. In this equation $X_{k,l}$ is a parameter in $Y_{k,l,m}$ though this will not be 
explicitly written.  Since we assume that the dynamics of the $y$ variables is ergodic,
we may choose $Y_{k,l,0} = Y_{k,l-1,m}$ for all $k$ and $l$. For convenience however, we will set $Y_{k,l,0}=0.$ As for the $X$ variables, our notation
implies $Y_{k,0,m}=Y_{k-1,L,m}$ for all $m.$

The
existence, under appropriate assumptions, of an invariant measure, $\mu_x^{\dt},$ of  the
numerical scheme \eqref{eq:microsolver}, follows from results in \cite{higham,roberts}.
The measure $\mu_x^{\dt}$ is an approximation to $\mu_{x}.$
This suggests estimating the function $\ba$ by 
\begin{equation}
\begin{aligned}
& A(X_{k,l}) =  
\frac{1}{M} \sum_{m=1}^M
a(X_{k,l},Y_{k,l,m}).
\end{aligned}
\label{eq:AandB}
\end{equation}
Since we will frequently encounter this form of trajectory averaging in the sequel, we define the
symbol
$$
\left[\mathcal{S}^M h\right] (y_1,\dots,y_m) = \frac{1}{M} \sum_{m=1}^M h(y_m),
$$
where $h$ is some function defined on $\mathbb{R}^{d_y}.$  We will henceforth omit from our notation,
the 
dependence of $\mathcal{S}^M h$ on the variables $(y_1,\dots,y_m).$
Equations \eqref{eq:macro}, \eqref{eq:microsolver}, and
\eqref{eq:AandB} define the multiscale integration scheme.
We note here that, in fact, since for our filtering application of the multiscale integration scheme we are
interested in reproducing the distribution of the process
\eqref{eq:system} and not actual trajectories, we can replace
the average in \eqref{eq:AandB} by evaluation at the single point
$Y_{k,l,M}.$  Since expression \eqref{eq:AandB} is only marginally more expensive and corresponds more naturally to the method of averaging, we will not bother with this simplification.

Suppose that the joint distribution of
$Y_{k,l,m}$ given that $X_{k,l} = x$ and $Y_{k,l,0}=0$ is 
$\widetilde Q^{\dt,M}_x(dy_1,\dots,dy_M).$
We can define the measure $\mu^M_x$ by,
$$
\mu_x^M h = \widetilde Q^{\dt,M} \mathcal{S}^M h(y_m)
= \frac{1}{M}\sum_{m=1}^M \mathbf{E}_x h(Y_{k,l,m}).
$$
where the subscript $x$ on the expectation emphasizes the dependence of the random variables $Y_{k,l,m}$
on the parameter $x.$

By using trajectory averages over the fast dynamics to approximate integrals over $\mu_x$
 of the observation density 
 as well as the coefficient $\bar a$ we can define a particle filter
 which approximates the Rao-Blackwellization step of the previous section. 
The following kernel is an approximation of expression \eqref{eq:kernelaveragedRB},
\begin{multline}\label{eq:kernelfinal}
\bar{H}_{k}^{\Dt,\dt,M}f\,(x_{s_{k-1}},z_k)
= \int \bar Q_{\Delta s}^{\Dt,\dt,M}(x_{s_{k-1}},dx_{s_k})\\
\times\mu_{x_{s_k}}^M(dy)\,f(x_{s_k},y)\, g(x_{s_k},y,z_k).
\end{multline}

The next particle filter, defined in the following algorithm, corresponds to 
\eqref{eq:kernelfinal} and differs from Algorithm  \ref{avefilterrb}
in that we evolve the particles according to $\bar Q^{\Dt,\dt,M}_{\Ds}$ instead of 
$\bar Q_{\Ds}$ and, in the update Step, instead averaging over the measure $\mu_x$ we 
average over a trajectory of $Y^j_{k,l,m}.$
\begin{algorithm}{Multiscale particle filter}\label{finalfilter}

\begin{enumerate}
\item Sample $n$ i.i.d. vectors $X_0^j \; j=1,\ldots,n$ from $Q_0(dx,dy).$ Set $k=1.$

\item For each $j=1,\ldots,n$  evolve $X'_{k-1,0}= X^j_{k-1}$ according to
Equations \ref{eq:macro}, \ref{eq:microsolver}, and \ref{eq:AandB} to generate
samples  ${X'_k}^j = {X'}_{k-1,L}^j$, i.e. sample from  $\bQ^{\Dt,\dt,M}_{\triangle s}(X_{k-1}^j,dx').$
We use ${Y'}^j_{k-1,l,m}$ to denote the fast variables
evolved according to Eq. (\ref{eq:microsolver}) with parameter ${X'}^j_{k-1,l}.$  To initialize the
${Y'}^j$ variables
at each time Step of Eq. (\ref{eq:macro}) we choose,
${Y'}^j_{k,l,0} = 0.$
\item   \label{item:returnDtM} For each $j=1,\ldots,n:$

    \begin{enumerate}
        \item \label{item:unresolved} Evolve ${Y'}_{k,0,0}^j=0$ according to (\ref{eq:microsolver}) with
        parameter ${X'_k}^j$ to get  ${Y'}_{k,0,m}^j\; \text{where} \; 1\leq m\leq M$.

        \item Calculate

    \begin{align*}
    \left[\mathcal{S}^Mfg\right]({X'_k}^j) &= \frac{1}{M}\sum_{m=1}^M f({X'_k}^j,Y_{k,0,m}^j)\\ 
	&\hspace{1cm}\times g({X'}_{k}^j,Y_{k,0,m}^j,z_k)\\
    \left[\mathcal{S}^Mg\right]({X'_k}^j) &=
 	\frac{1}{M}\sum_{m=1}^M \, g({X'_k}^j,Y_{k,0,m}^j,z_k).
    \end{align*}

                  \item Calculate
        \[
        \begin{split}
        W^n_k\,f &= \frac{1}{n}\sum_{j=1}^n \left[\mathcal{S}^Mfg\right]({X'_k}^j)\\
        W^n_k\,1 &= \frac{1}{n}\sum_{j=1}^n \left[\mathcal{S}^Mg\right]({X'_k}^j).
        \end{split}
        \]
       \end{enumerate}

\item   Define the probability measure $\Psi^n_k$ on $\Re^{d}$ with density
    \begin{multline*}
    \bar\Psi^{\Dt,\dt,M,n}_k(dx,dy) = \frac{1}{M n\,W^n_k\,1}\\
  	\times\sum_{j=1}^n \sum_{m=1}^M \,
    g({X'_k}^j,Y_{k,0,m}^j,z_k)\, \delta_{{X'_k}^j}(dx)\, \delta_{Y_{k,0,m}^j}(dy)
    \end{multline*}

\item \label{item:resample}If $k < N$: Set $ k \rightarrow k+1.$  Sample $n$ i.i.d. variables
$X_{k}^j$, $\; j=1,\ldots,n$ from $\Psi^n_k$, and  go to Step 2.
\end{enumerate}
We stop at the end of iteration $N$, and our approximation of $\Pi_k\,f$ is,
\[
\bar\Psi_k^{\Dt,\dt,M,n}\,f = \frac{W_k^n\,f}{W_k^n\,1} = \e_{\bar\Psi_k^{\Dt,\dt,M,n}} f
\]
\end{algorithm}

Since the $y$ components of the resampled particles after Step 4 are not used in Step 2, we
can, in practice, resample from the marginal density
\begin{equation}
 \label{altresample2}
\bar{\Psi}^{\Dt,\dt,M,n}_k(dx) = \frac{1}{n\,W^n_k\,1}
	\sum_{j=1}^n \,
     \left[\mathcal{S}^M g\right]\left({X'_k}^j\right)\, \delta_{{X'_k}^j}(dx).
\end{equation}
The two procedures give equivalent estimates and differ only in
that the work required for the resampling Step using \eqref{altresample2} will scale with $n$ instead of $nM$ (of course calculating the averaged weights requires $\mathcal{O}(nM)$ work).
In practice to initialize
${Y'}^j$ variables
at each time Step of Eq. (\ref{eq:macro}) we set,
${Y'}^j_{k,l,0} = {Y'}^j_{k,l-1,m}.$  This choice results in faster equilibration of the $Y'$ process.

While the details of the multiscale integration scheme do depend on the particular Markov process under study, algorithm \ref{finalfilter}
can be applied in the same form to a wide variety of multiscale Markov processes
(for example the system in section \ref{sec:jump}).

\section{Numerical results}\label{sec:numbers}

\subsection{A stochastic differential equation}\label{sec:sde}
\begin{figure}
 \centering
	\includegraphics[width=8cm]{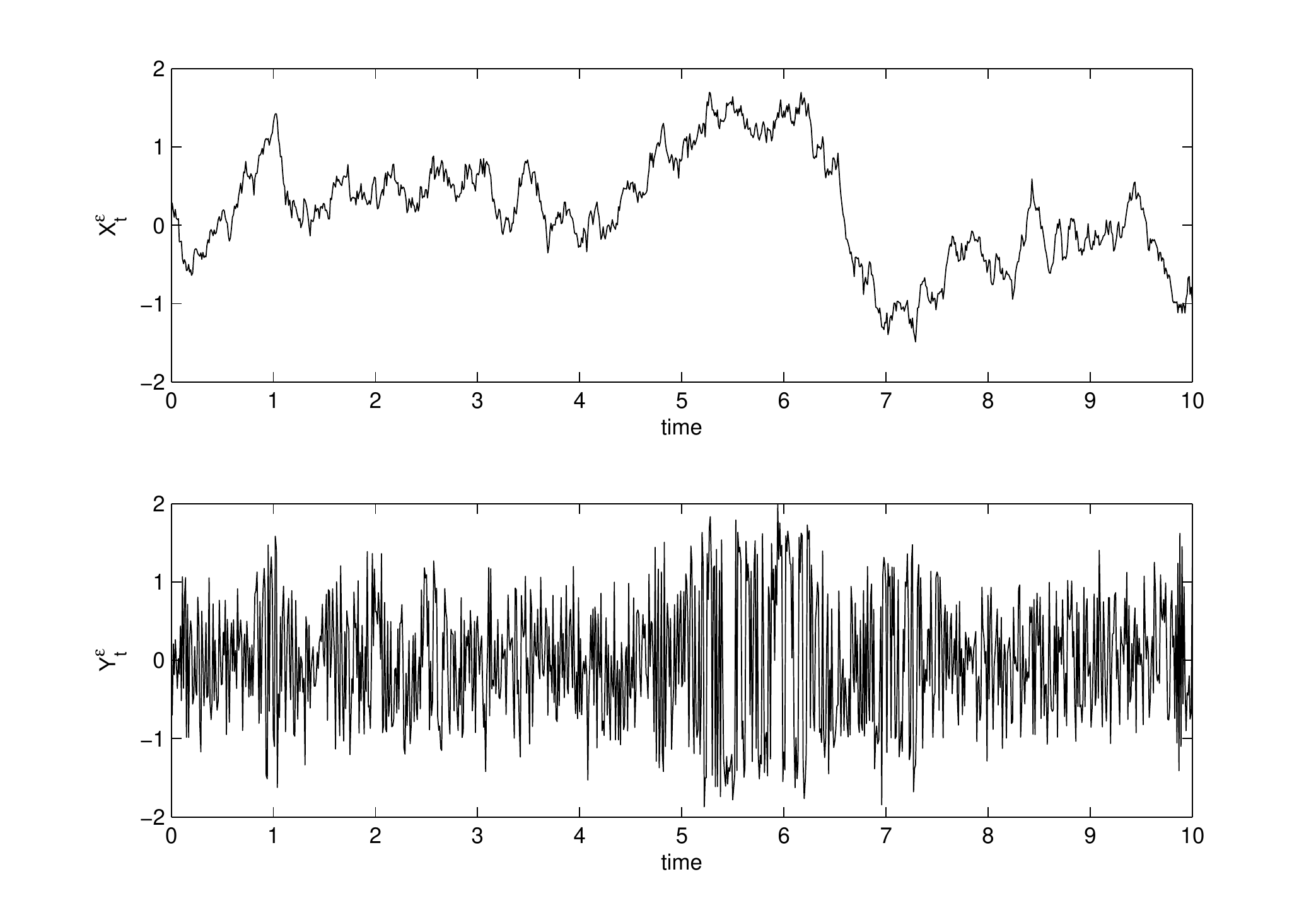}
  \caption{Trajectory of system \eqref{eq:numer}.}
  \label{fig:sdetraj}
  \end{figure}

\begin{figure}
 \centering
 \includegraphics[width=8cm]{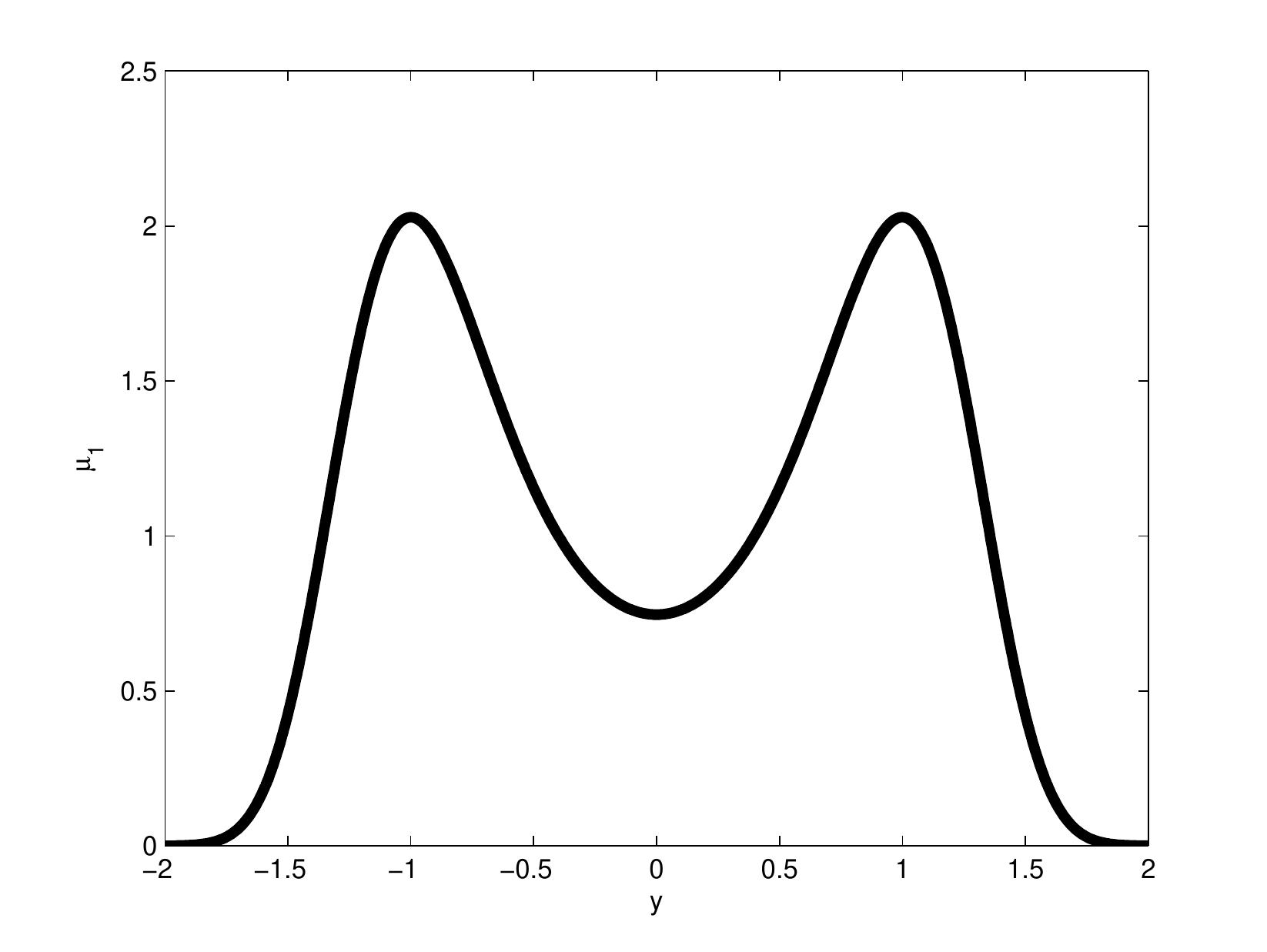}
  \caption{The density $\mu_x$ for $x=1.$}
  \label{fig:mu}
  \end{figure}

We present a simple numerical example which demonstrates the variance
reduction obtained through our algorithm. Consider the system given by,
\begin{equation}\label{eq:numer}
\begin{aligned}
dX^{\epsilon}_t  &= \left(Y^{\epsilon}_t - (X^{\epsilon}_t)^3
\right)
\,dt + dU_t &X^{\epsilon}_0 &\sim \mathcal{N}(0,1)\\
dY^{\epsilon}_t &= \frac{2}{\epsilon}\left((X^{\epsilon}_t)^2 - (Y^{\epsilon}_t)^2\right)Y^{\epsilon}_t\,dt +
\frac{1}{\sqrt{\eps}}\,dV_t  &Y^{\epsilon}_0 &\sim \mathcal{N}(0,1).
\end{aligned}
\end{equation}
The parameter $\eps$ in this example is set to $10^{-4}.$
A trajectory of the system is shown in Figure \ref{fig:sdetraj}.
The ergodic measure of the fast dynamics for this system is
known and has the bimodal density
$$
\mu_x(y) = \frac{e^{-(x^2-y^2)^2}}{\int e^{-(x^2-y^2)^2}dy}.
$$
A plot of $\mu_x$ for $x=1$ is shown in Figure \ref{fig:mu}.

The observations are given by,
\[
Z_k=Y^{\epsilon}_k+ \chi_k,
\]
where $\chi$ are independent Gaussian random variables with
mean 0 and standard deviation 0.1.  In the experiment below
we take as the realization of the observations, $z_k,$ the
trajectory shown in Figure \ref{fig:sdetraj} sampled at every
one unit of time.
We will compare the standard particle filtering algorithm with
Algorithm \ref{finalfilter}.  Both methods are run with
1000 particles. System \eqref{eq:numer} is discretized by the
Euler-Maruyama method with time step $\delta t=10^{-6}.$  The multiscale
integration scheme uses a time steps of size $\triangle t = 10^{-2}$
to evolve the reduced system \eqref{eq:macro} and of size $\delta t$
in the microscopic system \eqref{eq:microsolver}. 
With this choice of parameters Algorithm \ref{finalfilter}
runs in about half of the time of the standard particle filter.

\begin{figure}
 \centering
	\includegraphics[width=8cm]{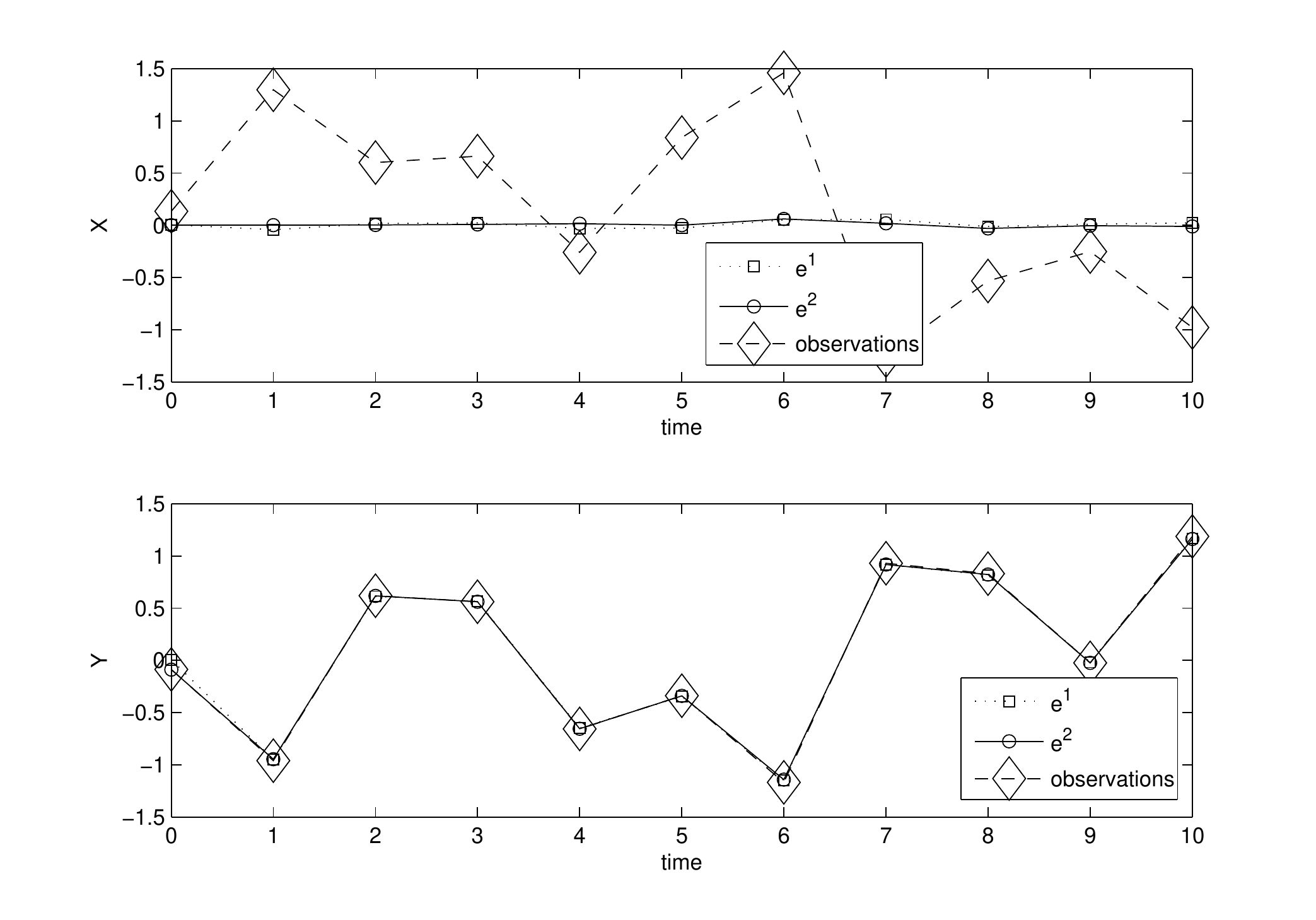}
  \caption{Trajectories of $e^1$ and $e^2.$}
  \label{fig:sdeest}
  \end{figure}

As in the definitions of the particle filters above, for any function $f$ define
$$
W_1f\,(k) = \sum_{j=1}^n f({{X_k^{\epsilon}}'}^j,{{Y^{\epsilon}}'}^j)
\,g({{X_k^{\epsilon}}'}^j,{{Y^{\epsilon}}'}^j,z_k)
$$
and
\begin{multline*}
W_2f\,(k) = \\ \frac{1}{n10^4}\sum_{j=1}^n\sum_{m=1}^{10^4} f({X_k'}^j,Y_{k,0,m}^j)\,g({X_k'}^j,Y_{k,0,m}^j,z_k).
\end{multline*}
In Figure \ref{fig:sdeest} we compare the pairs of estimators,
$$
e^1_x(k) = \frac{W_1 x\, (k)}{W_1 1\, (k)},\ \ \ \ 
e^1_y(k) = \frac{W_1 y\, (k)}{W_1 1\, (k)}
$$
and
$$
e^2_x(k) = \frac{W_2 x\, (k)}{W_2 1\, (k)},
\ \ \ \ 
e^2_y(k) = \frac{W_2 y\, (k)}{W_2 1\, (k)}
$$
Notice that the poor quality of the reconstruction of $X^{\epsilon}_t$ is not due to an error.  The symmetry
of the observation model and of the $Y^{\epsilon}_t$ dynamics,
implies that, 
in the limit $\epsilon\rightarrow 0,$ 
the true condition expectation of
$X^{\epsilon}_t,$ given any observations of $Y^{\epsilon}_t$ alone, will be 
identically zero.  Therefore, both estimators appear to be accurate.

In order  to compare the quality of the samples generated by the
two methods we compute the effective sample sizes of the
empirical measures produced by both methods,
$$
ess_1(k) = \frac{1000}{1+C^2_1(k)}
\quad\text{and}\quad ess_2(k) = \frac{1000}{1+C^2_2(k)}
$$
where
\begin{equation*}
C_1(k) =  \frac{1}{W_1 1\,(k)}
\sqrt{\frac{1}{n}\sum_{j=1}^n \left(g({{X_k^{\epsilon}}'}^j,{{Y^{\epsilon}}'}^j,z_k)
-W_1 1\,(k)\right)^2}
\end{equation*}
and
\begin{multline*}
C_2(k) = \frac{1}{W_2 1\,(k)}\sqrt{\frac{1}{n10^4}}\\
\times\sqrt{\sum_{j=1}^n\sum_{m=0}^{10^4}
\left(g({X_k'}^j_6(l),Y_{k,0,m}^j,z_k)
-W_2 1\,(k)\right)^2}.
\end{multline*}
The effective sample size is a common measure of the quality 
of a weighted empirical measure produced by an importance sampling scheme (see \cite{liu95}).  Very roughly, the effective sample size gives the number of independent samples from the target measure that would 
produce an unweighted empirical measure of equal quality.
When all of the weights are equal (i.e. the samples are drawn from the target measure itself) the effective sample size is
the total number of samples (in our case 1000).

\begin{figure}
 \centering
 \includegraphics[width=8cm]{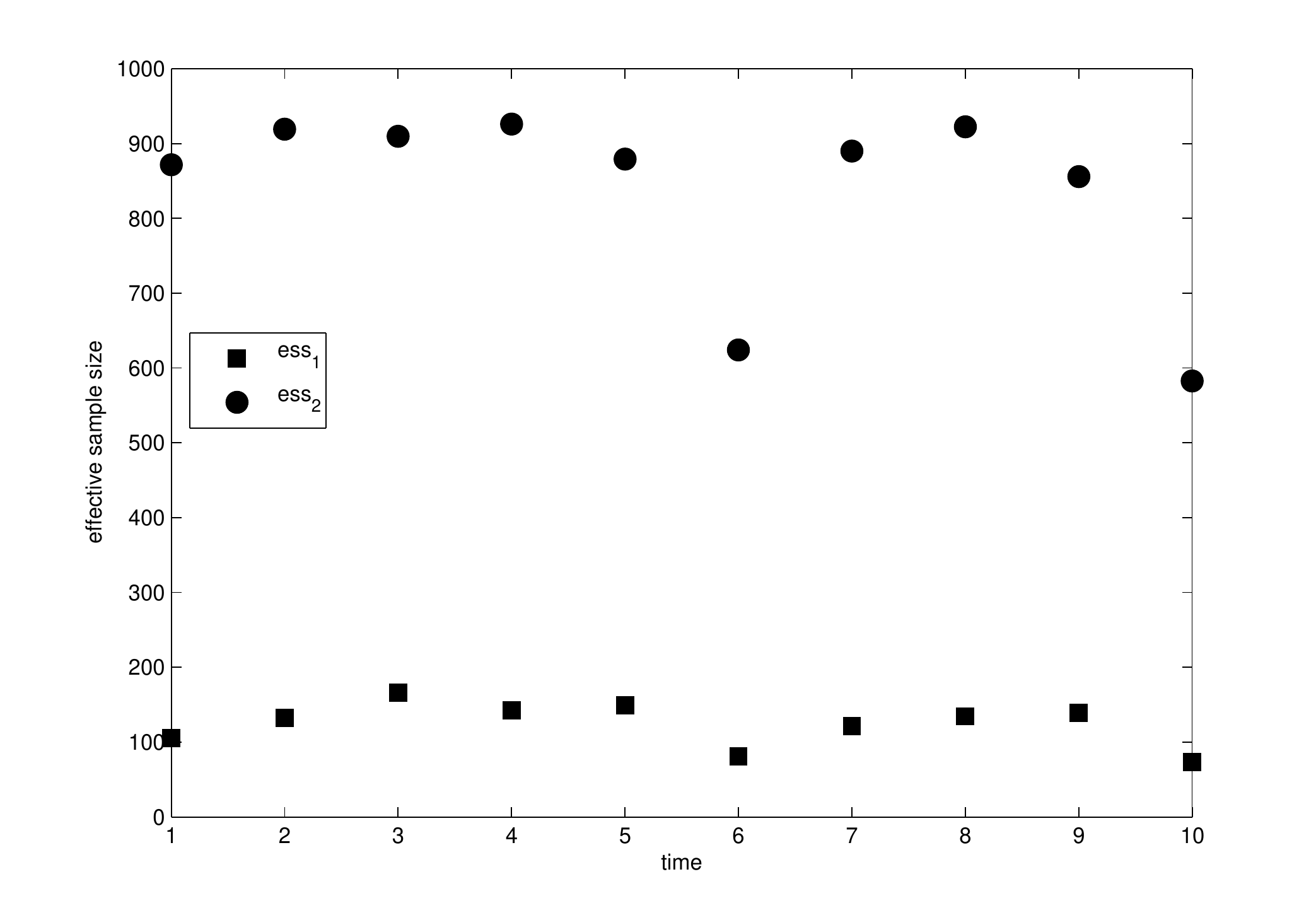}
  \caption{Trajectories of $ess_1$ and $ess_2.$}
  \label{fig:sdeess}
  \end{figure}

The trajectories of $ess_1$ and 
$ess_2$ are plotted in Figure \ref{fig:sdeess}.
As can be seen in the plot, the effective sample sizes generated by Algorithm \ref{finalfilter} are
nearly 10 times as large as those generated by  the standard particle filter.
This
indicates that there is some improvement in the quality of the empirical measure generated by Algorithm \ref{finalfilter}.

It is important to note while the cost of the two methods in this experiment are comparable,  in the limit $\eps\rightarrow 0$
a discretization of the system \eqref{eq:numer}
would require smaller and smaller time steps.  Thus the computational advantage for the multiscale
particle filter would become extreme in this limit.
The next example features a larger time scale separation 
and a correspondingly larger gain in efficiency due to multiscale
integration.

\subsection{A pure jump type Markov process}\label{sec:jump}

We now demonstrate our algorithm with a numerical example motivated by
chemical reactions. The stochastic dynamical behavior of a well
stirred mixture of $N$ molecular species that chemically interact
through $M$ reaction channels is accurately described by the
chemical master equation. The master equation is usually simulated
using the Stochastic Simulation Algorithm of Gillespie \cite{Gil76}.
 In
cellular systems where the small number of molecules of a few
reactant species sometimes necessitates a stochastic description of
the system's temporal behavior, chemical reactions often take place
on vastly different time scales. An exact stochastic simulation of
such a system will necessarily spend most of its time simulating the
more numerous fast reaction events.

For $N$ molecular species  the state of the system
$\bold{S}=\brk{S_1(t),\ldots,S_N(t)}$ is the number of  molecules
of each species
present at time $t$. The molecular populations $S_i,
i=1,\ldots,N$ are random variables. For each reaction channel $R_j,
j=1,\ldots,M$ we define the propensity function $a_j(\bold{S})$ and
the state change vector $v_j$. The propensity function is such that
$a_j(\bold{S})dt $ is the probability given $\bold{S}(t)=\bold{S}$
that one $R_j$ reaction will occur in the next infinitesimal time
interval $\Brk{t,t+dt}$. The state change vector $v_j$ is the change
in the number of $S_i$ molecules produced by one $R_j$ reaction.

The pathwise evolution law for the master equation is a jump type
Markov process on the non negative $N$-dimensional integer lattice
given by
\begin{multline}\label{eq:systemjump}
\left(
\begin{array}{c}
dS_1 \\
\vdots \\
dS_N
\end{array} \right)
=
 \left(
\begin{array}{ccc}
v_{1,1} & \ldots & v_{1,M} \\
\vdots &  & \vdots \\
v_{N,1} & \ldots & v_{N,M}
\end{array} \right)\\
\times
\left(
\begin{array}{c}
dP_1(a_1(\bold{S})) \\
\vdots \\
dP_M(a_M(\bold{S}))
\end{array} \right),
\end{multline}
where $P_j$ is a Poisson process with state dependent intensity
parameter $a_j(\bold{S})$.

\begin{figure}
 \centering
 \includegraphics[width=8cm]{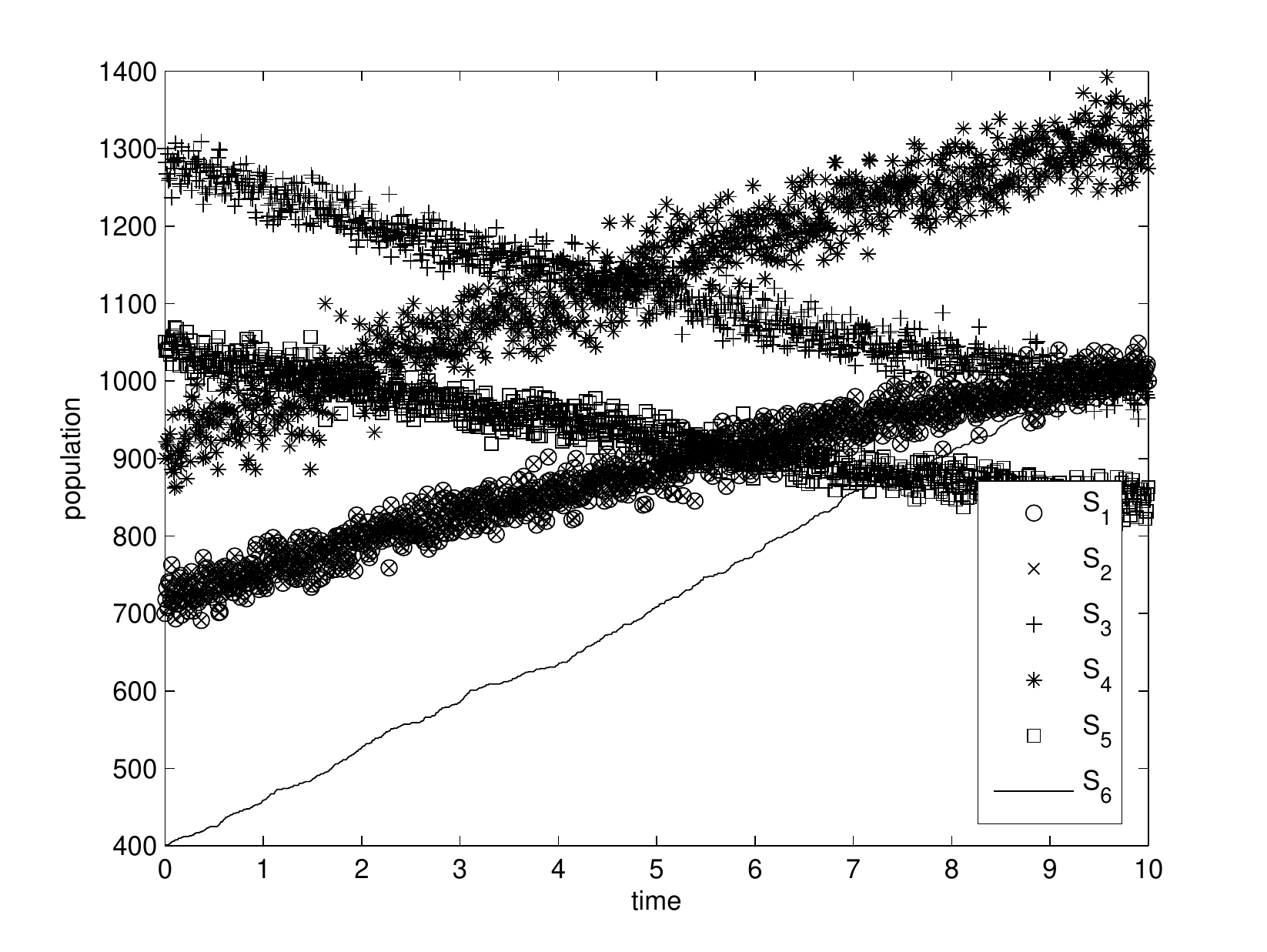}
  \caption{A trajectory of the jump process with initial condition
 (700,700,1300,900,1050,400).}
  \label{fig:jumptraj}
  \end{figure}

In our system we choose $N=6$ species and $M=5$ reaction channels.
Variables $S_1,\dots,S_5$ are the fast variables and $S_6$ is
the slow variable.
For the evolution of the fast variables we choose a simple fast biomolecular
(reversible) reaction
\begin{align*}
 S_1 + S_2 &\to_{ k_1} S_3\\
 S_3 +S_6 &\to_{k_2} S_1 +S_2 +S_6,
\end{align*}
and a fast (reversible) dimerization
\begin{align*}
 S_4 + S_4 &\to_{ k_3} S_5\\
 S_5 +S_6 &\to_{k_4} S_4 +S_4+S_6.
\end{align*}
The propensity functions are given by
\[
\begin{split}
\left(
\begin{array}{c}
a_1(\bold{S}) \\
a_2(\bold{S}) \\
a_3(\bold{S})\\
 a_4(\bold{S})
\end{array} \right)
&= \left(
\begin{array}{c}
k_1 S_1 S_2 \\
k_2 S_3 S_6\\
k_3 S_4 (S_4-1)\\
 k_4 S_5 S_6
\end{array} \right)
\end{split}
\]
where the reaction rates $k_1,\dots,k_5$ are specified below.
We will use the shorthand
$$
S_f = \left(S_1,\dots,S_5\right).
$$
For the slow variable $S_6$ we choose an external source
(spontaneous creation)
\[
S_3+S_4 \to_{k_5} S_3+S_4+S_6,
\]
which is coupled to the fast variables through   the slow reaction's
propensity function
\[
a_5(\bold{S})=k_5 S_3 S_4.
\] 
The state change vectors for the system just described are given 
by the matrix,
$$
\left(
\begin{array}{ccc}
v_{1,1} & \ldots & v_{1,5} \\
\vdots &  & \vdots \\
v_{6,1} & \ldots & v_{6,5}
\end{array} \right)
= \left(
\begin{array}{ccc}
-1,\ 1,\ 0,\ 0,\ 0\\
-1,\ 1,\ 0,\ 0,\ 0\\
\ 1,-1,\ 0,\ 0,\ 0\\
\ 0,\ 0,-2,\ 2,\ 0\\
\ 0,\ 0,\ 1,-1,\ 0\\
\ 0,\ 0,\ 0,\ 0,\ 1
\end{array}\right)
$$
Finally the reaction constants vector is given by
\[
\brk{k_1, \ldots,k_5} = \brk{1000, 1000, 1000, 1000, 5\times 10^{-5}}
\]
A trajectory of this system is plotted in Figure \ref{fig:jumptraj}.
One can clearly see the separation in time scale between the slow and 
fast variables.  By examination of the propensities, $a_i,$ above one can see that the $S_f$ variables
evolve on a time scale of roughly $10^{-9}$ while $S_6$ evolves on a time scale of
roughly $10^{-2}.$

The initial populations, $S_1(0),\dots,S_6(0)$ are drawn according to
\begin{equation*}
\left(S_1(0),\dots,S_6(0)\right) = (700,700,1300,900,1050,400) + \eta
\end{equation*}
where $\eta$ is a vector of independent uniformly distributed random 
variables in the interval $[-10,10].$
The observations are taken every 1 unit of time from time 0 to time 10 
and are modeled by
$$
Z(l) = S(l) + \chi(l)
$$
where the $\chi(l)$ are independent vectors of independent, mean 0, 
standard deviation 5, Gaussian random variables.  While this noise model is
somewhat arbitrary, it can be considered to model measurement error.
In our experiments we take $z(l)$ to be equal to the
particular trajectory shown in Figure \ref{fig:jumptraj} at time
$l$ for $l=0,\dots,10.$

Because the standard particle filter is too expensive to run (almost
1000 times the cost of the multiscale particle filter) it is not available
for comparison.  However, we can still investigate the variance reduction
aspect of the method by comparing two filters that both use the multiscale 
integration scheme during the prediction step, but that handle the update 
step in different ways.  These two procedures are described in the next paragraph.  For the details of the implementation of
the multiscale integration scheme applied here please see \cite{CGP05,elv07}.

The two particle filters applied here correspond roughly to algorithms
\ref{avefilter} and \ref{avefilterrb}.  In the first, at iteration
$l$, we apply
the multiscale integration scheme during the prediction step (thereby 
generating an 
approximate sample, ${S'}^j_6(l)$ from $\bar Q_{\Ds}$) 
then we choose one sample ${S'}^j_f$ approximately from the measure 
$\mu_{{S'}^j_6(l)}$ and proceed as in the rest of algorithm
 \ref{avefilter}.  The latter sampling step is accomplished by sampling the
end point 
from a long ($10^{-4}$ time units) equilibrated trajectory of the fast dynamics with the
value ${S'}^j_6(l)$ set as a parameter.  The second method
(corresponding to \ref{avefilterrb}) proceeds in exactly the
same way except that in the second step we sample $10^4$ points ${S'}^j_{f,m},$
at equal time intervals,
from a long (again $10^{-4}$ time units) equilibrated trajectory of the fast dynamics with the
value ${S'}^j_6(l)$ set as a parameter.  The points ${S'}^j_{f,m}$
are then used to approximate the two averages appearing in Step 3 of
algorithm \ref{avefilterrb}.
With the severe scale separation present in this problem, the variance of
the first method just described (no likelihood weight averaging) 
is an accurate representation of the variance
of the standard particle filter.  The increased cost of the averaging
in the second method is negligible.

Both particle filters are tested with 1000 particles.  
The resulting estimators,
$$
e^1_i(l) = \frac{W_1 s_i\,(l)}{W_1 1\, (l)}\qquad
\text{and}\qquad e^2_i(l) = \frac{W_2 s_i\,(l)}{W_2 1\, (l)}
$$
where, for any function $f,$
$$
W_1f\,(l) = \frac{1}{n}\sum_{j=1}^n f({S'}^j_6(l),{S'}^j_f)\,g({S'}^j_6(l),{S'}^j_f,z_l),
$$
and
\begin{multline*}
W_2f\,(l) = \\ \frac{1}{n10^4} \sum_{j=1}^n\sum_{m=1}^{10^4} f({S'}^j_6(l),{S'}^j_{f,m})\,g({S'}^j_6(l),{S'}^j_{f,m},z_l).
\end{multline*}
\begin{figure} \centerline{
{\bf a.}
\includegraphics[width=8cm]{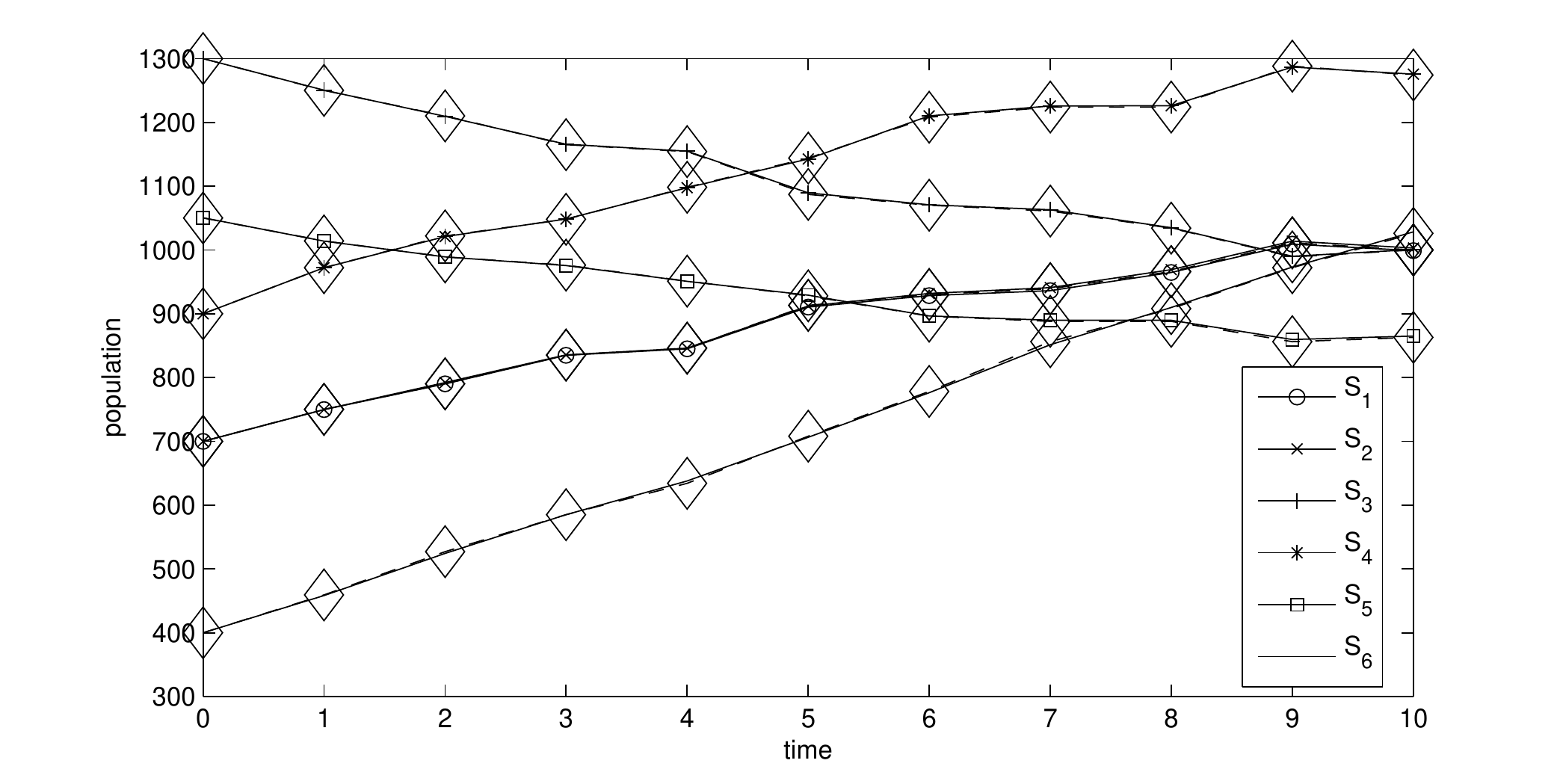}}
\centerline{
{\bf b.}
\includegraphics[width=8cm]{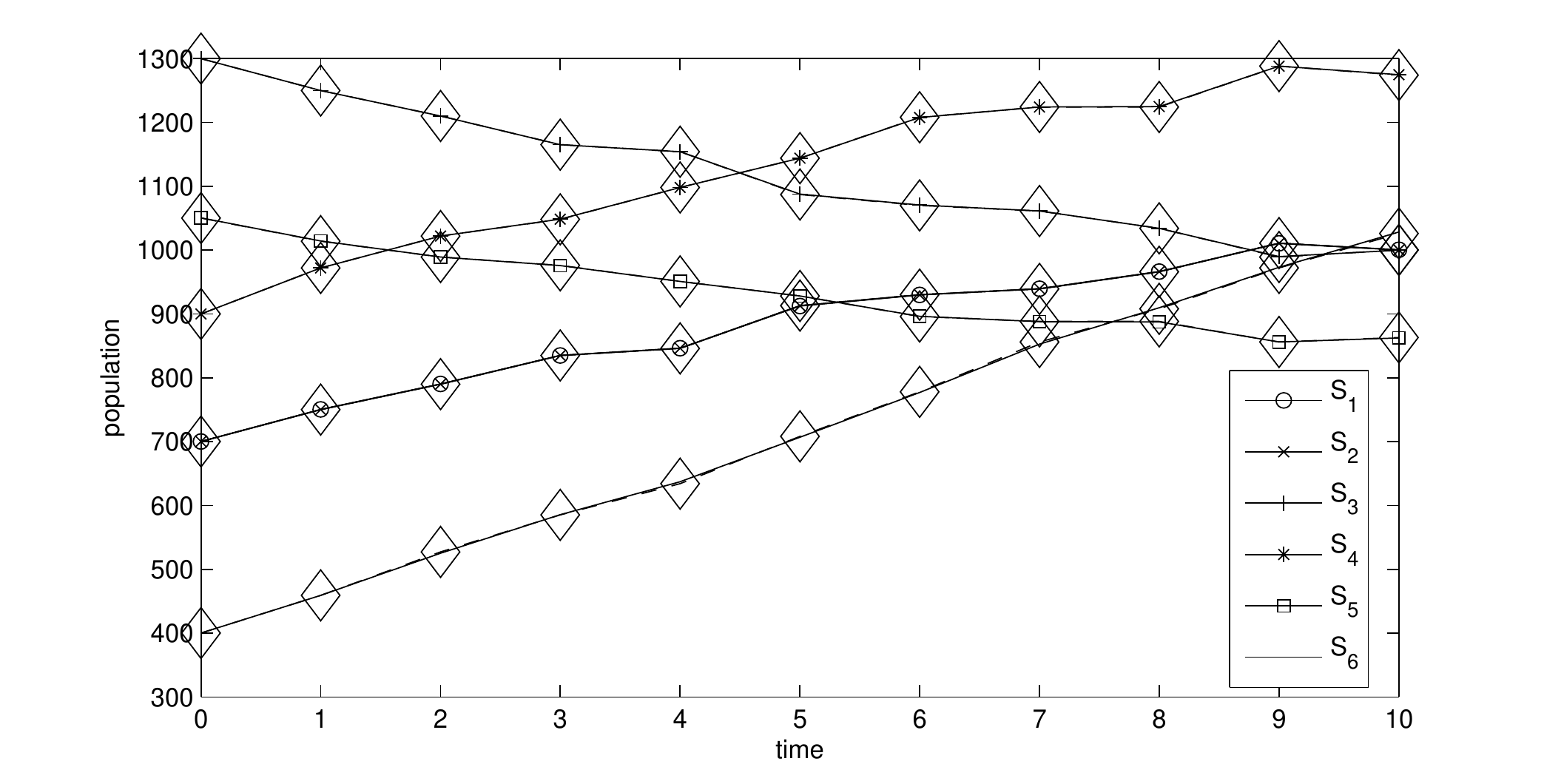}}
\caption{{\bf a.}   Trajectory of estimate $e^1$.  
{\bf b.} Trajectory of estimate $e^2$.  In both figures
the observations are overlaid with diamonds connected by
dashed lines.}
\label{fig:jumpest}
\end{figure}
of $S_i$ are plotted in Figure \ref{fig:jumpest} along with the hidden
trajectory of $S_6$.
Clearly both methods accurately reconstruct the path of $S.$
Upon close inspection one can see that the estimate $e^2$ is 
slightly closer to the actual trajectory (which is the same as the
observations in this example).  However Figure \ref{fig:jumpest} is by no means
conclusive evidence that that $e^2$ provides a better estimate.

\begin{figure}
 \centering
 \includegraphics[width=8cm]{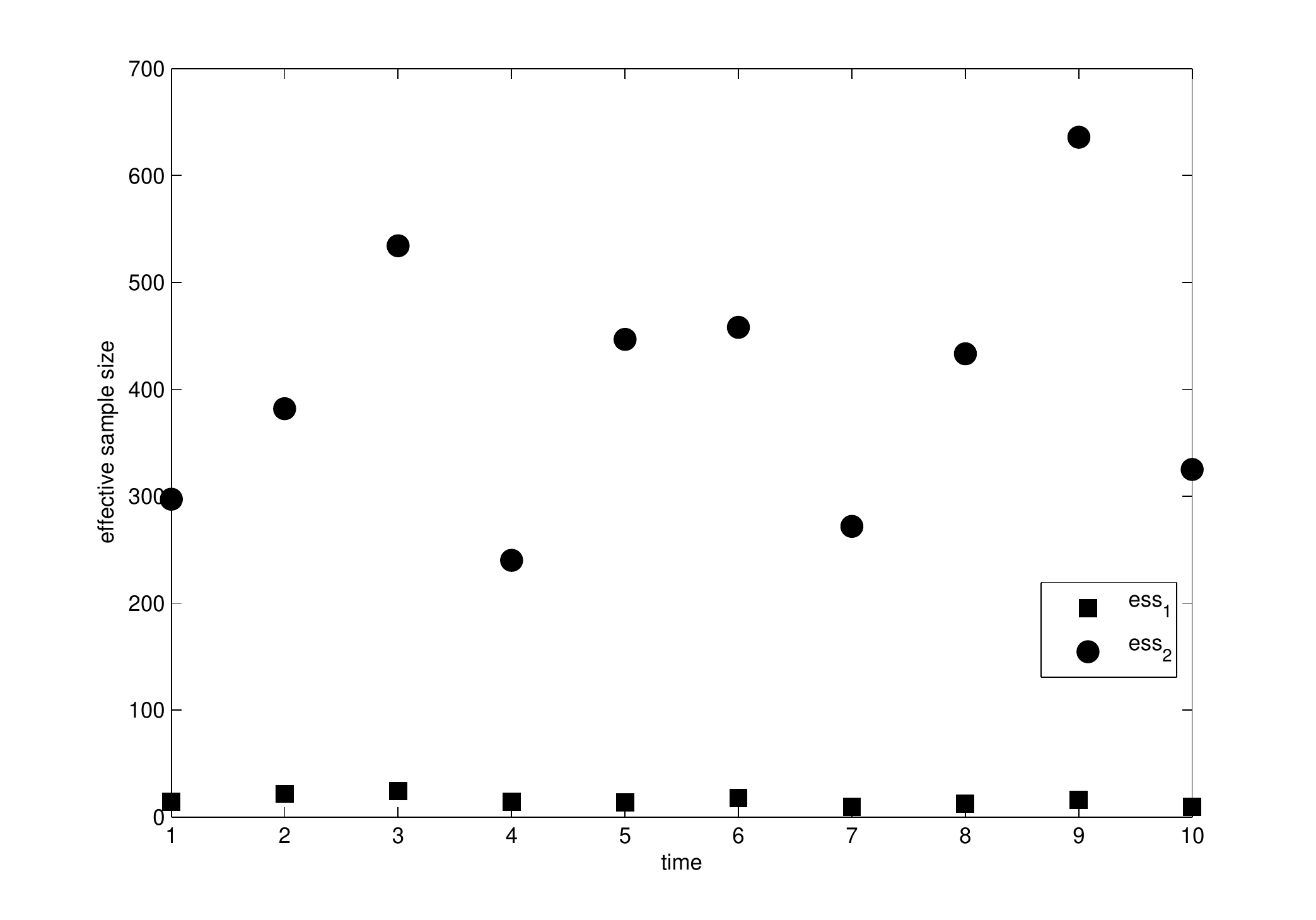}
  \caption{Trajectories of $ess_1$ and $ess_2.$}
  \label{fig:jumpess}
  \end{figure}

The gain in efficiency of the empirical measure generated
using the averaged weights does become clear when we compare the effective sample sizes of the
empirical measures produced by both methods,
$$
ess_1 = \frac{1000}{1+C^2_1}
\quad\text{and}\quad ess_2 = \frac{1000}{1+C^2_2}
$$
where
\begin{equation*}
C_1(l) = \frac{1}{W_11\,(l)}\sqrt{\frac{1}{n}\sum_{j=1}^n \left(g({S'}^j(l),{S'}_f^j,z_l)
-W_11\,(l)\right)^2}
\end{equation*}
and
\begin{multline*}
C_2(l) = \frac{1}{W_21\,(l)}\sqrt{\frac{1}{n10^4}}\\
\times\sqrt{\sum_{j=1}^n\sum_{m=0}^{10^4}
\left(g({S'}^j_6(l),{S'}^j_{f,m},z_l)
-W_21\,(l)\right)^2}.
\end{multline*}
The trajectories of $ess_1$ and 
$ess_2$ are plotted in Figure \ref{fig:jumpess}.
As can be seen in the plot, the effective sample sizes generated by the particle filter with averaging are as much as 100 times greater than those generated by  the particle filter without the averaging step.
This
indicates that the improvement in the quality of the empirical measure generated by the
particle filter with likelihood weight averaging is significant.
It is also clear that, due to the large time scale separation in this system, any filtering method that does not explicitly address this in the prediction stage of the filter (by, for example, multiscale integration) will be extremely slow.

\section{Conclusions}\label{sec:conclusions}
We have presented an algorithm that combines dimensional reduction
and approximate Rao-Blackwellization to create an efficient reduced variance
particle filter  for systems that exhibit time scale separation. We
tested the algorithm on two systems with large time scale
separations and the results are encouraging.

Our method does not require that any of the distributions involved are nearly 
Gaussian or degenerate.  Furthermore, the cost of our method does not increase as the time
scale separation is increased.  We are not aware of any competitive alternative
with these features.


With minor modifications
our particle filter can be applied in a slightly more general setting. 
Occasionally one is interested in multiscale systems for which it is impossible to
explicitly fix the slow variable during evolution.  For example, one may not
know the laws governing the evolution of the system but can generate
sample trajectories (by laboratory experiment).
In these cases
our particle filter remains valid.
Roughly, the fact that 
the
 slow variable does not deviate much on the time scale of the fast
variables allows one to replace evolution of the fast dynamics by evolution
of the full system (see \cite{kevrekidis03}).  Note also, that the variance reduction
technique discussed here applies to any importance sampling problem for a multiscale
Markov process.


\section{Acknowledgments} We are grateful to Professor A.J. Chorin for
suggesting the use of dimensional reduction in the context of filtering. 
We would like to thank Dr. Y. Shvets
and the anonymous referees for helpful suggestions regarding this
manuscript.   This work was supported in part by the National
Science Foundation under Grant DMS 04-32710, and by the Director,
Office of Science, Computational and Technology Research, U.S.\
Department of Energy under Contract No.\ DE-AC03-76SF000098.

\bibliographystyle{unsrt}
\bibliography{GSW07.bib}

\end{document}